\documentclass[graybox]{svmult}

\usepackage{type1cm}        
\usepackage{makeidx}         
\usepackage{graphicx}        
\usepackage{multicol}        
\usepackage[bottom]{footmisc}
\usepackage{makecell}

\usepackage{newtxtext}       %
\usepackage[varvw]{newtxmath}       

\makeindex             

\usepackage[utf8]{inputenc}
\usepackage{amsmath}
\usepackage{amsfonts}
\usepackage{mathrsfs}
\usepackage{bm}
\usepackage{times}
\usepackage{soul}
\usepackage{graphicx}
\usepackage{subcaption}
\usepackage{color, colortbl}
\definecolor{Gray}{gray}{0.95}
\newcommand{\comment}[1]{}

\def\bf{\bm{f}}

\def\div{\text{\rm div}}

\def\hrot{H(\text{\rm rot};\Omega)}

\def\V0{V_0^{\rm edge}(E)} 

\begin{document}

\title*{Coarse spaces for virtual element methods on irregular 3D subdomain decompositions}
\titlerunning{Coarse spaces for VEM on irregular 3D subdomain decompositions}
\author{Ana Aguilar-Pineda, Luis F. Amey, Adrián Angulo-Paniagua and Juan G. Calvo}
\institute{Juan G. Calvo \at CIMPA - Escuela de Matemática, Universidad de Costa Rica, Costa Rica\\
\email{juan.calvo@ucr.ac.cr}
\and Ana Aguilar-Pineda \at Escuela de Matemática, Universidad de Costa Rica, Costa Rica\\
\email{ana.aguilarpineda@ucr.ac.cr}
\and Luis F. Amey \at  Escuela de Matemática, Universidad de Costa Rica, Costa Rica\\
\email{luis.amey@ucr.ac.cr}
\and Adrián Angulo-Paniagua \at  Escuela de Matemática, Universidad de Costa Rica, Costa Rica\\
\email{adrian.angulopaniagua@ucr.ac.cr}}
%
%
\maketitle

\abstract*{We present a two-level overlapping Schwarz preconditioner for three-dimensional problems discretized with the Virtual Element Method. Our approach handles general polyhedral meshes and irregular subdomains, extending the applicability of previous methods. Numerical experiments show robust performance with respect to the number of subdomains and mesh parameters, with condition-number bound comparable to classical finite element results. While alternative methods such as FETI-DP and BDDC are available, the simplicity and competitiveness of overlapping additive Schwarz methods underscore the practical significance of our contribution.}

\section{Introduction}
\label{sec:1}

Given a bounded polyhedral domain $\Omega\subset \mathbb{R}^3$, we seek $u\in H_0^1(\Omega)$ such that
\begin{equation}\label{weakForm}
a(u,v)\ := \int_{\Omega} \rho \nabla u\cdot \nabla v = \int_{\Omega} f v\quad \forall v\in H_0^1(\Omega),
\end{equation} 
where $\rho \in L^\infty(\Omega)$ is a positive function that is uniformly bounded from below. We recall that $H_0^1(\Omega)$ is the space of functions in $L^2(\Omega)$ with gradient also in $L^2(\Omega)$ and vanishing trace. 

The discretization of the variational problem by finite element or virtual element methods (FEM or VEM, respectively) generally leads to large and ill-conditioned linear systems, especially when the coefficient $\rho$ exhibits strong variations or when nonuniform meshes are employed. In this setting, the design of efficient preconditioners is essential to ensure the scalability and rapid convergence of the iterative methods used to solve such systems. Among the available strategies, additive Schwarz preconditioners constitute a particularly appealing family because they are built from independent local solvers defined on each subdomain, which simplifies their implementation and makes them suitable for complex geometries. However, achieving scalability requires the inclusion of an appropriate global coarse space, whose construction must be carefully designed to ensure robustness with respect to mesh refinement and coefficient variations. Compared with more sophisticated approaches such as BDDC or FETI-DP, additive Schwarz schemes offer a favorable balance between computational performance and practical simplicity, while retaining solid theoretical properties and good parallelization capabilities.

Early studies on overlapping Schwarz methods focused on tetrahedral or cubic subdomains; see, e.g., \cite[Chapter 3]{T&W}. The study of irregular subdomains for FEM in 2D was initiated for problems posed in $H^1(\Omega)$ in \cite{D&K&W_DDLessRegSubd,MR2436090,MR2436089,Widlund_2008}. Such theory has been studied in 2D for John and Jones domains and the rates of convergence of the preconditioners are determined by parameters that characterize such domains. Coarse spaces consist of vertex and edge functions that are extended as discrete harmonic functions into the interiors of the subdomains. 

Authors in \cite{D&W_AltCoarseSpace} defined coarse functions at each subdomain vertex by assigning specific values on the interface and extending them harmonically into the interior for 2D. These ideas were later generalized to VEM in \cite{MR3815657,MR3913656} for 2D problems, where projectors replaced the harmonic extensions. Such coarse spaces form a partition of unity, which can be used to construct robust domain decomposition methods for high-contrast multiscale problems on irregular domains by solving suitable eigenvalue problems to enrich the coarse space \cite{CALVO2024112909}. Moreover, a small coarse space is defined even for irregular subdomains in \cite{MR3686806}, but numerical examples for 3D problems were reported only for cubic subdomains. For problems posed in $\hrot$ in 2D and $H(\div;\Omega)$ in 3D, see \cite{MR3407231,MR3454359} and \cite{MR3033012, MR3243012}, respectively. 

In 3D, however, to the best of the authors' knowledge, no results exist for the systematic construction of global coarse functions on general irregular partitions with appropiate energy estimates for VEM. Relevant works for FEM include \cite{NadeemJimack}, which studies a parallel implementation of a two-level additive Schwarz preconditioner for tetrahedral elements; \cite{UReading}, which analyzes two-level Schwarz methods for matrices from p-version FEM on triangular and tetrahedral meshes with the coarse level given by the lowest-order finite element space; and \cite{improving3D}, which introduces an energy-minimizing coarse space of reduced size, with basis functions forming a partition of unity on subdomain interfaces and extended as discrete harmonic functions into the interior with cubic subdomains.

In this paper, we present a two-level overlapping Schwarz preconditioner for problem \eqref{weakForm} discretized with VEM in three dimensions. To the best of our knowledge, there are no theoretical results for additive preconditioners for the linear system that arises from \eqref{weakForm} when VEM are used and irregular subdomains are considered. Our method allows us to handle irregular subdomains and general polygonal meshes, and applies to a broader range of material properties and subdomain geometries than previous studies. Similarly as in previous studies, we assume that $\rho$ is constant on each subdomain of the decomposition.

A theoretical bound for the condition number $\kappa$ of a two-level overlapping Schwarz preconditioner for FEM, based on two reduced versions of the GDSW coarse space in three dimensions \cite{MR1302680}, has the form
\begin{equation} \label{eq:kappa}
\kappa\leq C \left(1+\frac{H}{\delta}\right) \left(1+\log\frac{H}{h}\right)^\alpha,
\end{equation}
where $\alpha \in \{0,1,2\}$ and $C$ is a positive constant independent of the number of subdomains; see \cite{improving3D} for further details, where authors present results for cubic subdomains.

We numerically observe a bound comparable to \eqref{eq:kappa} for our VEM preconditioner. We remark that there are different DDM such as FETI-DP and BDDC methods; see \cite{MR3735828,MR4101370} for 2D and 3D studies related to our problem. Nevertheless, the simplicity of implementing an {overlapping additive} Schwarz algorithm with competitive results gives relevance to our work.

The rest of this paper is organized as follows. In Section \ref{sec:vem}, we briefly describe the VEM for our model problem \eqref{weakForm}. Section \ref{sec:algorithm} presents the two-level overlapping additive Schwarz method and the construction of our coarse space. Numerical results are given in Section \ref{sec:numericalEx}. Finally, Section \ref{sec:conc} contains our conclusions and closing remarks.

\section{The virtual element method} \label{sec:vem}
For the sake of completeness, we briefly describe the 3D VEM; we follow \cite{MR2997471,beirao}. We remark some numerical implementations for 2D \cite{Sutton2016,MR4544725} and 3D problems \cite{mVEM,VEMcomp}. Let $\mathcal{T}_h$ be a decomposition of $\Omega$ into general polyhedral elements $E$; see \cite{BeiraoStab} for assumptions and study on the stability of the method. In Figure \ref{fig:1}, we present a hexagonal prism and a Voronoi mesh. For simplicity, we consider the lowest order virtual space. We denote by $\mathbb{P}_k(\mathcal{O})$ the set of polynomials of degree at most $k$ defined in $\mathcal{O}$ (where $\mathcal{O}$ could represent an edge, face or element), and denote by $\{m_{\bm{\alpha}}\}$ the set of scaled monomials (for edges, faces or elements).

First, we define the virtual space $W(f)$ for a face $f\in\partial E$. A function $u\in W(f)$ satisfies: (i) $u|_e \in \mathbb{P}_1(e)$ for each edge $e\in\partial f$; (ii) $u|_{\partial f}\in C^0(\partial f)$; (iii) $\Delta u \in \mathbb{P}_1(f)$, and (iv) $\int_f u m_\alpha = \int_f \Pi_f^\nabla u m_\alpha$ for a base of linear polynomials defined on $f$. The projector $\Pi_f^\nabla u$ is defined as the $L^2-$projection of the gradient of $u$, where $u$ and $\Pi_f^\nabla u$ have the same nodal average.

A virtual element function $u$ belongs to the local virtual space $V(E)$ if: (i) $v|_e\in\mathbb{P}_1(e)\ \forall \text{ edge } e\in \partial E$, (ii) $\ v|_f \in W(f) \ \forall \text{ face } f \in \partial E$, (iii) $v|_{\partial E}\in C^0(\partial E)$, and (iv) $\Delta v=0$. It is clear that $\mathbb{P}_1(E) \subseteq V(E)$. The degrees of freedom are chosen as the values of $u$ at the vertices of $E$.

The projector $\Pi^\nabla : V(E)\to \mathbb{P}_1(E)$ is defined by (i) $\int_E \nabla p \cdot \nabla \Pi_E^\nabla u=\int_E \nabla p \cdot \nabla u$ for all $p\in\mathbb{P}_1(E)$, and (ii) $\sum u(\bm{x}_i) = \sum \Pi_E^\nabla u(\bm{x}_i)$, where the sum goes over the vertices $\bm{x}_i$ of $E$. The global conforming virtual element space is then
$$ V_h := \{ v \in H^1_0(\Omega) \ : \ v|_E \in V(E) \ \ \forall E \in \mathcal{T}_h \}.$$

To discretize the weak formulation of the Poisson equation, the discrete bilinear form is then defined similarly as the 2D case:
$$
a_h^E(u,v) := (\nabla \Pi_E^\nabla u_h, \nabla \Pi_E^\nabla v_h)_E \;+\; S^E\big((I-\Pi_E^\nabla)u_h,(I-\Pi_E^\nabla)v_h\big),
$$
where $S^E(\cdot,\cdot)$ is a stabilization bilinear form computable from the nodal values of $u_h$ and $v_h$, typically chosen as $S^E(u_h,v_h)=h_E \sum u_h(\bm{x}_i) v_h(\bm{x}_i)$, where the sum goes over all the vertices $\bm{x}_i$ of $E$. The global matrix is then assembled via 
$$ a_h(u_h,v_h) := \sum_{E \in \mathcal{T}_h} a_h^E(u_h,v_h),$$
together with a discrete load term obtained via the $L^2$-projection of $f$ onto $\mathbb{P}_1(E)$.  
This yields a conforming discretization of the Poisson problem on general polyhedral meshes with optimal approximation properties. We omit implementation details; for the 3D discretization we have used the library \texttt{mVEM}; see \cite{mVEM}.

\section{Coarse space} \label{sec:algorithm}
\begin{figure}[t!]
\sidecaption[t]
\includegraphics[scale=.395]{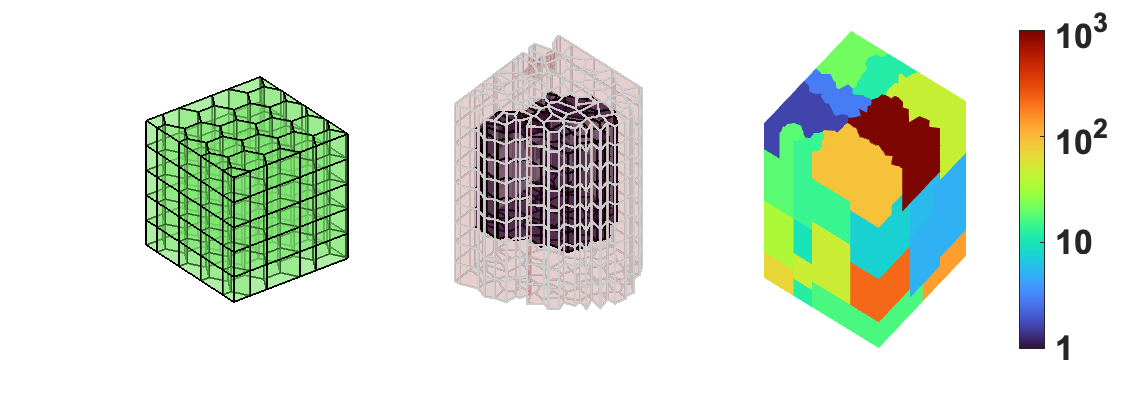}
\caption{(left) A polyhedral mesh, (middle) a subdomain $\Omega_i$ and its overlapping subdomain $\Omega_i'$, and (right) a piecewise coefficient $\rho$ with $\rho|_{\Omega_i} = \rho_i \in [1,10^3]$.}
\label{fig:1}
\end{figure}

We consider the standard two-level additive Schwarz preconditioner; see \cite[Chapter 3]{T&W} for further details. Consider a coarse partition $\mathcal{T}_H$ consisting of subdomains $\{\Omega_i\}_{i=1}^N$, where each subdomain is a union of elements; see Figure \ref{fig:1}. The interface of the decomposition $\{\Omega_i\}_{i=1}^N$ is defined as the union of the local boundaries $\partial \Omega_i$, excluding the external boundary $\partial \Omega$. Equivalence classes can be associated with the subdomain faces, edges, and vertices.

\begin{figure}[t!]
\sidecaption[t]
\includegraphics[scale=.398]{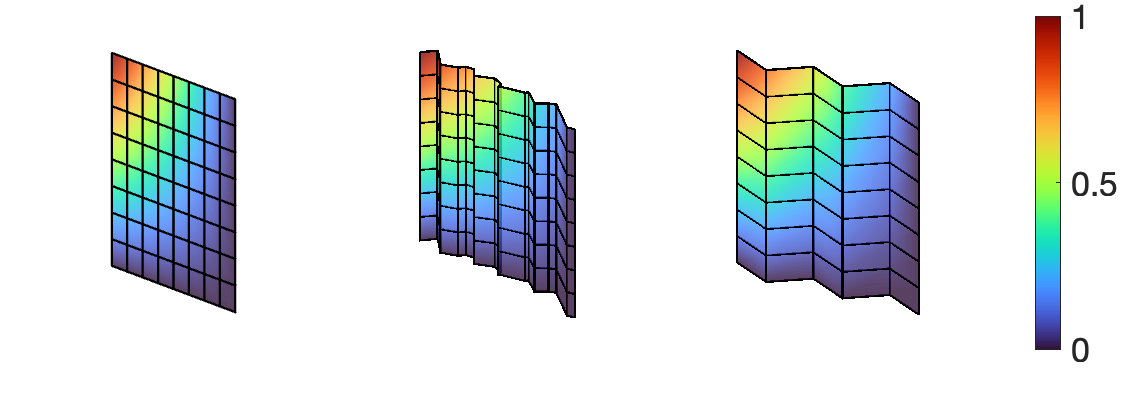}
\caption{Subdomain faces for (left) square, (middle) Voronoi and (right) hexagonal meshes with subdomains based on the incenter of the elements. $R_0^T \varphi_{\bm{x}_0}$ is shown for each subdomain face.}
\label{fig:coarseFnA}
\end{figure}

\begin{figure}[t!]
\sidecaption[t]
\includegraphics[scale=.398]{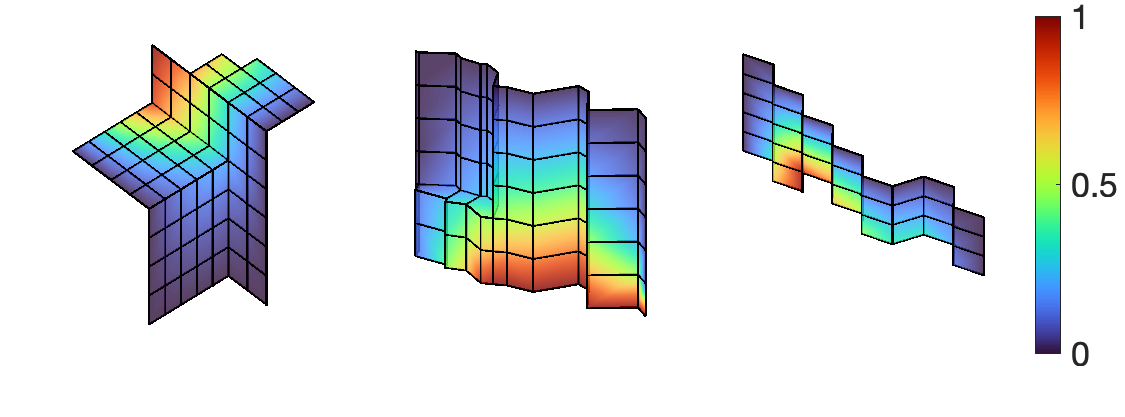}
\caption{Subdomain faces for (left) square, (middle) Voronoi and (right) hexagonal meshes with METIS subdomains. $R_0^T \varphi_{\bm{x}_0}$ is shown for each subdomain face.}
\label{fig:coarseFnB}
\end{figure}

A subdomain face $\mathcal{F}_{ij}$ corresponds to the degrees of freedom associated with nodes lying in the interior of the intersection of the boundaries of two neighboring subdomains $i$ and $j$, excluding any edges on the boundary of the face. If this intersection consists of multiple disjoint components, each component is considered a separate face; see Figures \ref{fig:coarseFnA} and \ref{fig:coarseFnB} for typical subdomain faces.

A subdomain edge $\mathcal{E}_{ij,kl}$ is defined as the interior of the intersection of the closure of faces $\mathcal{F}_{ij}$ and $\mathcal{F}_{kl}$, typically consisting of nodes shared by three or more subdomains. Subdomain vertices are the endpoints of all subdomain edges, excluding the vertices in $\partial \Omega$, since homogeneous Dirichlet boundary conditions are assumed.

The coarse space $V_0$ could be defined as the virtual element space associated to the coarse mesh $\mathcal{T}_H$ composed by subdomains, but its dimension would in general be too large. To obtain a coarse space of moderate dimension, we instead define a reduced space as follows. For each subdomain vertex $\bm{x}_0$, we define a coarse function $\varphi_{\bm{x}_0}$ that extends ideas presented in \cite{D&W_AltCoarseSpace} for 2D problems first and then generalized to 2D VEM in \cite{MR3815657,MR3913656}. This function satisfies $\varphi_{\bm{x}_0}(\bm{x}_i)=\delta_{i0}$ for all subdomain vertices $\bm{x}_i$, and vanishes at nodes on any subdomain edge that does not include $\bm{x}_0$ as an endpoint. If $\mathcal{E}$ is an edge with endpoints $\bm{x}_0$ and $\bm{x}_1$, then $\varphi_{\bm{x}_0}$ varies linearly along the direction from $\bm{x}_0$ to $\bm{x}_1$; see \cite{MR3913656}.

For nodes on a subdomain face $\mathcal{F}$ with $\bm{x}_0$ on its boundary, $\varphi_{\bm{x}_0}$ is defined as the solution of the discrete Laplace equation $\Delta_\mathcal{F} \varphi_{\bm{x}_0}=0$, subject to the boundary values on $\partial \mathcal{F}$; see Figures \ref{fig:coarseFnA} and \ref{fig:coarseFnB}. For this discretization, we use the 2D VEM on each face. The face projectors are required for the 3D stiffness matrix, and hence are computed only once. In this construction, we obtain a partition of unity on the interface of the decomposition. 

Finally, the function is extended to the interior of the subdomains via discrete harmonic extension. Each function $\varphi_{\bm{x}_0}$ corresponds to a column of the extension operator $R_0^T$, which maps a function from the reduced coarse space to $V_h$. Following \cite{MR3815657}, it could be possible to consider polynomial approximations to the interior of the subdomains in order to obtain faster approximations.

We then construct overlapping subdomains $\Omega_i' \supset \Omega_i$ by adding layers of elements that are external to $\Omega_i$, and we will denote by $\delta_i$ the minimum width of the region $\Omega_i'\setminus\Omega_i$. For each subdomain $i$, $1\leq i\leq N$, we define the local virtual space
\begin{equation*} 
V_i:= \left\lbrace v\in H^1_0(\Omega_i'):  v|_{E}\in V(E)\ \forall\ E\subset \Omega_i'\right\rbrace.
\end{equation*}
The degrees of freedom of a function $v_i \in V_i$ are given by its values at the fine nodes lying in the interior of $\Omega_i'$. We also introduce the natural operators $R_i^T : V_i \rightarrow V_h$, $1 \le i \le N$, defined by zero extension from the subdomain $\Omega_i'$ to the global domain $\Omega$. We remark that the matrix representation of these operators $R_i^T$ are not required to be assembled; we only require the indices of internal nodes of $\Omega_i'$.

We finally consider the two-level additive overlapping Schwarz preconditioner
\begin{equation} \label{additiveOp}
P_{ad}:=\sum_{i=0}^N P_i = A^{-1}_{ad} A,\ \text{with}\ A^{-1}_{ad} = \sum_{i=0}^N R_i^T (R_iAR_i^T)^{-1} R_i,
\end{equation}
where we consider exact solvers for each subspace for simplicity; see \cite[Chap. 2]{T&W}.

\section{Numerical results} \label{sec:numericalEx}
We present numerical results for the two-level additive overlapping Schwarz preconditioner \eqref{additiveOp}. We solve the resulting linear systems using the preconditioned conjugate gradient method to a relative residual tolerance of $10^{-6}$. We estimate the condition number $\kappa(P_{ad})$ and compute the number of iterations $I$ for each experiment; see results in Table \ref{tab:results1}. We present numerical experiments for dependence on $N$ (scalability), $H/\delta$ and $H/h$, for cubic, Voronoi and hexagonal-prism elements. We include two sets of coefficients: (i) $\rho=1$ in the whole domain, and (ii) $\rho=\rho_D$, a piecewise constant function with $\rho|_{\Omega_i}= \rho_i\in [1,10^3]$. Numerical experiments were performed in \texttt{MATLAB} on a MacBook Pro equipped with an Apple M3 Pro chip and 18 GB of RAM, and the VEM discretizacion was obtained with \texttt{mVEM} \cite{mVEM}.

We  confirm the linear growth in the condition number as we increase $H/\delta$ and  we observe no significant dependence on the parameters $N$ and $H/h$, as expected from previous bounds for the condition number of the preconditioner system as \eqref{eq:kappa}. We remark that for the case of cubic elements, our method recovers the usual FEM space in the coarse mesh. 

\begin{table}[t]
\centering
\caption{Number of iterations $I$ and condition number $\kappa$ (in parentheses) for cubic, Voronoi, and hexagonal–prism triangulations with $N$ subdomains. Three sets of tests are reported, in which $N$, $H/\delta$, and $H/h$ are varied independently while keeping the other two parameters fixed. Both $\rho=1$ and $\rho=\rho_D$ (a discontinuous, piecewise-constant $\rho$ on each subdomain) are considered. The quantity $|V_0|$ denotes the dimension of the coarse space.}
\label{tab:results1}
\begin{tabular}{l|c|c|c|c|c|c|c|c|c|}
\cline{2-10}
 & \multicolumn{3}{c|}{Cubes} 
 & \multicolumn{3}{c|}{Voronoi} 
 & \multicolumn{3}{c|}{Hexagonal prisms} \\ \cline{2-10}

 & \multicolumn{1}{c|}{\makecell{$\rho = 1$ \\ $I(\kappa)$}}
 & \multicolumn{1}{c|}{\makecell{$\rho = \rho_D$ \\ $I(\kappa)$}}
 & \multicolumn{1}{c|}{$|V_0|$}

 & \multicolumn{1}{c|}{\makecell{$\rho = 1$ \\ $I(\kappa)$}}
 & \multicolumn{1}{c|}{\makecell{$\rho = \rho_D$ \\ $I(\kappa)$}}
 & \multicolumn{1}{c|}{$|V_0|$}

 & \multicolumn{1}{c|}{\makecell{$\rho = 1$ \\ $I(\kappa)$}}
 & \multicolumn{1}{c|}{\makecell{$\rho = \rho_D$ \\ $I(\kappa)$}}
 & $|V_0|$ \\ \hline

\rowcolor[HTML]{EFEFEF}
\multicolumn{1}{|r|}{$N$} 
 & \multicolumn{9}{l|}{\cellcolor[HTML]{EFEFEF}Test 1: $H/h = 8$, $H/\delta = 2$} \\ \hline

\multicolumn{1}{|r|}{$4^3$}
 & 18(10.7) & 23(10.6) & 27
 & 27(13.3) & 23(13.6) & 54
 & 25(12.4) & 26(12.2) & 54 \\ \hline

\multicolumn{1}{|r|}{$6^3$}
 & 18(11.9) & 26(12.0) & 125
 & 26(13.3) & 23(13.2) & 250
 & 25(11.6) & 28(12.0) & 250 \\ \hline

\multicolumn{1}{|r|}{$8^3$}
 & 18(13.1) & 27(13.0) & 343
 & 28(14.3) & 24(17.3) & 686
 & 25(12.7) & 30(12.8) & 686 \\ \hline

\multicolumn{1}{|r|}{$10^3$}
 & 18(14.8) & 28(13.8) & 729
 & 30(14.8) & 35(16.3) & 1458
 & 25(13.0) & 30(13.3) & 1458 \\ \hline

\rowcolor[HTML]{EFEFEF}
\multicolumn{1}{|r|}{$H/\delta$}
 & \multicolumn{9}{l|}{\cellcolor[HTML]{EFEFEF}Test 2: $N=4^3$, $H/h = 12$} \\ \hline

\multicolumn{1}{|r|}{3}
 & 19(13.1) & 24(13.6) & 27
 & 23(10.6) & 20(12.3) & 54
 & 22(10.7) & 26(12.8) & 54 \\ \hline

\multicolumn{1}{|r|}{4}
 & 19(15.3) & 25(15.7) & 27
 & 24(11.6) & 21(18.4) & 54
 & 23(11.2) & 28(13.6) & 54 \\ \hline

\multicolumn{1}{|r|}{6}
 & 20(20.4) & 27(18.7) & 27
 & 28(15.5) & 23(18.4) & 54
 & 26(14.6) & 31(15.7) & 54 \\ \hline

\multicolumn{1}{|r|}{12}
 & 22(34.3) & 33(31.6) & 27
 & 39(32.2) & 34(37.1) & 54
 & 33(27.5) & 38(22.0) & 54 \\ \hline

\rowcolor[HTML]{EFEFEF}
\multicolumn{1}{|r|}{$H/h$}
 & \multicolumn{9}{l|}{\cellcolor[HTML]{EFEFEF}Test 3: $N=3^3$, $H/\delta = 4$} \\ \hline

\multicolumn{1}{|r|}{8}
 & 17(14.4) & 22(13.6) & 8
 & 24(11.2) & 26(10.9) & 16
 & 22(11.9) & 23(10.4) & 16 \\ \hline

\multicolumn{1}{|r|}{16}
 & 18(13.1) & 22(13.4) & 8
 & 25(10.7) & 26(10.5) & 16
 & 23(11.1) & 25(11.3) & 16 \\ \hline

\multicolumn{1}{|r|}{24}
 & 18(11.2) & 23(12.3) & 8
 & 27(11.0) & 26(11.2) & 20
 & 23(10.8) & 27(11.7) & 16 \\ \hline

\end{tabular}
\end{table}

\section{Conclusions} \label{sec:conc}

Although a theoretical analysis of the proposed method is still pending, the numerical results demonstrate its robustness with respect to the number of subdomains and the ratio $H/h$, while exhibiting approximately linear growth with the relative overlap. As a direction for future work, alternative approaches to harmonic extensions inside the subdomains could be explored, for instance by using projections intrinsic to the VEM spaces, which may lead to significant computational savings, as in \cite{MR3913656}. Furthermore, in the case of high-contrast or multiscale coefficients $\rho$, the partition of unity provided by the coarse space is expected to enable the use of strategies similar to those in \cite{CALVO2024112909}, potentially yielding robust and efficient solvers in such challenging scenarios.

\begin{acknowledgement}
The authors gratefully acknowledge the institutional support for project C5054 provided by the \textit{Vicerrectoría de Investigación, Universidad de Costa Rica}, through the Research Center in Pure and Applied Mathematics (CIMPA) and the School of Mathematics.
\end{acknowledgement}

\bibliographystyle{siam}
\bibliography{biblio}

@article {MR4101370,
    AUTHOR = {Bertoluzza, Silvia and Pennacchio, Micol and Prada, Daniele},
     TITLE = {F{ETI}-{DP} for the three dimensional virtual element method},
   JOURNAL = {SIAM J. Numer. Anal.},
  FJOURNAL = {SIAM Journal on Numerical Analysis},
    VOLUME = {58},
      YEAR = {2020},
    NUMBER = {3},
     PAGES = {1556--1591},
      ISSN = {0036-1429,1095-7170},
   MRCLASS = {65N30 (65N55)},
  MRNUMBER = {4101370},
MRREVIEWER = {Feng\ Wang},
       DOI = {10.1137/18M1233303},
       URL = {https://doi.org/10.1137/18M1233303},
}

@article {MR3735828,
    AUTHOR = {Bertoluzza, Silvia and Pennacchio, Micol and Prada, Daniele},
     TITLE = {B{DDC} and {FETI}-{DP} for the virtual element method},
   JOURNAL = {Calcolo},
  FJOURNAL = {Calcolo. A Quarterly on Numerical Analysis and Theory of
              Computation},
    VOLUME = {54},
      YEAR = {2017},
    NUMBER = {4},
     PAGES = {1565--1593},
      ISSN = {0008-0624,1126-5434},
   MRCLASS = {65N30 (65N55)},
  MRNUMBER = {3735828},
MRREVIEWER = {Christos\ Kravvaritis},
       DOI = {10.1007/s10092-017-0242-3},
       URL = {https://doi.org/10.1007/s10092-017-0242-3},
}

@article {MR1302680,
    AUTHOR = {Dryja, Maksymilian and Smith, Barry F. and Widlund, Olof B.},
     TITLE = {Schwarz analysis of iterative substructuring algorithms for
              elliptic problems in three dimensions},
   JOURNAL = {SIAM J. Numer. Anal.},
  FJOURNAL = {SIAM Journal on Numerical Analysis},
    VOLUME = {31},
      YEAR = {1994},
    NUMBER = {6},
     PAGES = {1662--1694},
      ISSN = {0036-1429},
   MRCLASS = {65N55 (65N22 65N30)},
  MRNUMBER = {1302680},
MRREVIEWER = {B.\ Kellogg},
       DOI = {10.1137/0731086},
       URL = {https://doi.org/10.1137/0731086},
}

@incollection {MR3243012,
    AUTHOR = {Oh, Duk-Soon},
     TITLE = {An alternative coarse space method for overlapping {S}chwarz
              preconditioners for {R}aviart-{T}homas vector fields},
 BOOKTITLE = {Domain Decomposition Methods in Science and Engineering {XX}},
    SERIES = {Lect. Notes Comput. Sci. Eng.},
    VOLUME = {91},
     PAGES = {361--367},
 PUBLISHER = {Springer, Heidelberg},
      YEAR = {2013},
      ISBN = {978-3-642-35274-4; 978-3-642-35275-1},
   MRCLASS = {65N30 (65N55)},
  MRNUMBER = {3243012},
       DOI = {10.1007/978-3-642-35275-1\_42},
       URL = {https://doi.org/10.1007/978-3-642-35275-1_42},
}

@article {MR3033012,
    AUTHOR = {Oh, Duk-Soon},
     TITLE = {An overlapping {S}chwarz algorithm for {R}aviart-{T}homas
              vector fields with discontinuous coefficients},
   JOURNAL = {SIAM J. Numer. Anal.},
  FJOURNAL = {SIAM Journal on Numerical Analysis},
    VOLUME = {51},
      YEAR = {2013},
    NUMBER = {1},
     PAGES = {297--321},
      ISSN = {0036-1429,1095-7170},
   MRCLASS = {65N30 (65F08 65N55)},
  MRNUMBER = {3033012},
MRREVIEWER = {Huazhong\ Tang},
       DOI = {10.1137/110838868},
       URL = {https://doi.org/10.1137/110838868},
}

@article {MR3686806,
    AUTHOR = {Dohrmann, Clark R. and Widlund, Olof B.},
     TITLE = {On the design of small coarse spaces for domain decomposition
              algorithms},
   JOURNAL = {SIAM J. Sci. Comput.},
  FJOURNAL = {SIAM Journal on Scientific Computing},
    VOLUME = {39},
      YEAR = {2017},
    NUMBER = {4},
     PAGES = {A1466--A1488},
      ISSN = {1064-8275,1095-7197},
   MRCLASS = {65N55 (65F08 65F10 65N30)},
  MRNUMBER = {3686806},
MRREVIEWER = {Benjamin\ Wi-Lian\ Ong},
       DOI = {10.1137/17M1114272},
       URL = {https://doi.org/10.1137/17M1114272},
}

@article{CALVO2024112909,
title = {Robust domain decomposition methods for high-contrast multiscale problems on irregular domains with virtual element discretizations},
journal = {J. Comput. Phys.},
volume = {505},
pages = {112909},
year = {2024},
issn = {0021-9991},
doi = {https://doi.org/10.1016/j.jcp.2024.112909},
url = {https://www.sciencedirect.com/science/article/pii/S002199912400158X},
author = {Juan G. Calvo and Juan Galvis}
}

@InProceedings{improving3D,
author="Heinlein, Alexander
and Klawonn, Axel
and Rheinbach, Oliver
and Widlund, Olof B.",
editor="Bj{\o}rstad, Petter E.
and Brenner, Susanne C.
and Halpern, Lawrence
and Kim, Hyea Hyun
and Kornhuber, Ralf
and Rahman, Talal
and Widlund, Olof B.",
title="Improving the Parallel Performance of Overlapping {S}chwarz Methods by Using a Smaller Energy Minimizing Coarse Space",
booktitle="Domain Decomposition Methods in Science and Engineering XXIV ",
year="2018",
publisher="Springer International Publishing",
address="Cham",
pages="383--392",
isbn="978-3-319-93873-8"
}

@incollection {MR2436089,
    AUTHOR = {Dohrmann, Clark R. and Klawonn, Axel and Widlund, Olof B.},
     TITLE = {A family of energy minimizing coarse spaces for overlapping
              {S}chwarz preconditioners},
 BOOKTITLE = {Domain Decomposition Methods in Science and Engineering
              {XVII}},
    SERIES = {Lect. Notes Comput. Sci. Eng.},
    VOLUME = {60},
     PAGES = {247--254},
 PUBLISHER = {Springer, Berlin},
      YEAR = {2008},
      ISBN = {978-3-540-75198-4},
   MRCLASS = {65N55 (65N22)},
  MRNUMBER = {2436089},
       DOI = {10.1007/978-3-540-75199-1\_28},
       URL = {https://doi.org/10.1007/978-3-540-75199-1_28},
}

@incollection {MR2436090,
    AUTHOR = {Dohrmann, Clark R. and Klawonn, Axel and Widlund, Olof B.},
     TITLE = {Extending theory for domain decomposition algorithms to
              irregular subdomains},
 BOOKTITLE = {Domain Decomposition Methods in Science and Engineering
              {XVII}},
    SERIES = {Lect. Notes Comput. Sci. Eng.},
    VOLUME = {60},
     PAGES = {255--261},
 PUBLISHER = {Springer, Berlin},
      YEAR = {2008},
      ISBN = {978-3-540-75198-4},
   MRCLASS = {65N55},
  MRNUMBER = {2436090},
       DOI = {10.1007/978-3-540-75199-1\_29},
       URL = {https://doi.org/10.1007/978-3-540-75199-1_29},
}

@article {MR2997471,
    AUTHOR = {Beir\~ao da Veiga, L. and Brezzi, F. and Cangiani, A. and
              Manzini, G. and Marini, L. D. and Russo, A.},
     TITLE = {Basic principles of virtual element methods},
   JOURNAL = {Math. Models Methods Appl. Sci.},
  FJOURNAL = {Mathematical Models and Methods in Applied Sciences},
    VOLUME = {23},
      YEAR = {2013},
    NUMBER = {1},
     PAGES = {199--214},
      ISSN = {0218-2025},
   MRCLASS = {65N06},
  MRNUMBER = {2997471},
MRREVIEWER = {Bo\AA !`ko S. Jovanovi\"A},
       DOI = {10.1142/S0218202512500492},
       URL = {http://dx.doi.org/10.1142/S0218202512500492},
}

@article {MR4544725,
    AUTHOR = {Herrera, C\'esar and Corrales-Barquero, Ricardo and
              Arroyo-Esquivel, Jorge and Calvo, Juan G.},
     TITLE = {A numerical implementation for the high-order 2{D} virtual
              element method in {MATLAB}},
   JOURNAL = {Numer. Algorithms},
  FJOURNAL = {Numerical Algorithms},
    VOLUME = {92},
      YEAR = {2023},
    NUMBER = {3},
     PAGES = {1707--1721},
      ISSN = {1017-1398,1572-9265},
   MRCLASS = {65Y15 (65M60)},
  MRNUMBER = {4544725},
       DOI = {10.1007/s11075-022-01361-4},
       URL = {https://doi.org/10.1007/s11075-022-01361-4},
}

@Article{Sutton2016,
author="Sutton, Oliver J.",
title="{The virtual element method in 50 lines of MATLAB}",
JOURNAL = {Numer. Algorithms},
year="2017",
month="August",
day="01",
volume="75",
number="4",
pages="1141--1159",
}

@article {beirao,
    AUTHOR = {Beir\~ao da Veiga, L. and Brezzi, F. and Marini, L. D. and
              Russo, A.},
     TITLE = {The hitchhiker's guide to the virtual element method},
   JOURNAL = {Math. Models Methods Appl. Sci.},
  FJOURNAL = {Mathematical Models and Methods in Applied Sciences},
    VOLUME = {24},
      YEAR = {2014},
    NUMBER = {8},
     PAGES = {1541--1573},
      ISSN = {0218-2025},
   MRCLASS = {65N30},
  MRNUMBER = {3200242},
MRREVIEWER = {Nasser H. Sweilam},
       DOI = {10.1142/S021820251440003X},
       URL = {http://dx.doi.org/10.1142/S021820251440003X},
}

@article {MR3454359,
    AUTHOR = {Calvo, Juan G.},
     TITLE = {A {BDDC} algorithm with deluxe scaling for {$H({\rm curl})$}
              in two dimensions with irregular subdomains},
   JOURNAL = {Math. Comp.},
  FJOURNAL = {Mathematics of Computation},
    VOLUME = {85},
      YEAR = {2016},
    NUMBER = {299},
     PAGES = {1085--1111},
      ISSN = {0025-5718},
   MRCLASS = {65N30 (35Q60 65F10 65N55)},
  MRNUMBER = {3454359},
MRREVIEWER = {Marius Ghergu},
       DOI = {10.1090/mcom/3028},
       URL = {https://doi.org/10.1090/mcom/3028},
}

@article {MR3407231,
    AUTHOR = {Calvo, Juan G.},
     TITLE = {A two-level overlapping {S}chwarz method for {$H(\rm curl)$}
              in two dimensions with irregular subdomains},
   JOURNAL = {Electron. Trans. Numer. Anal.},
  FJOURNAL = {Electronic Transactions on Numerical Analysis},
    VOLUME = {44},
      YEAR = {2015},
     PAGES = {497--521},
   MRCLASS = {65N30 (35Q60 65N55)},
  MRNUMBER = {3407231},
}

@Inproceedings{Widlund_2008,
TITLE="Accommodating Irregular Subdomains in Domain Decomposition Theory",
AUTHOR="O.B.~Widlund",
booktitle={Domain Decomposition Methods in Science and Engineering XVIII},
series={Lecture Notes in Computational Science and Engineering},
publisher={Springer-Verlag},
PAGES="87--98",
volume="70",
YEAR="2009",
editor={Bercovier, M. and Gander, M. J. and Kornhuber, R. and O.B.~Widlund}
}

@ARTICLE{D&K&W_DDLessRegSubd,
   AUTHOR = {Dohrmann, Clark R. and Klawonn, Axel and Widlund, Olof B.},
   title = "{Domain decomposition for less regular subdomains: Overlapping Schwarz in two dimensions}",
   journal = "SIAM J. Numer. Anal.",
   volume = 46,
   pages = "2153-2168",
   year = 2008
   }

@ARTICLE{D&W_AltCoarseSpace,
   AUTHOR = {Dohrmann, Clark R. and Widlund, Olof B.},
   title = "{An alternative coarse space for irregular subdomains and an overlapping Schwarz algorithm for scalar elliptic problems in the plane}",
   journal = "SIAM J. Numer. Anal.",
   volume = 50,
   pages = "2522-2537",
   year = 2012
   }

@BOOK{T&W,
   author = "A.~Toselli and O.B.~Widlund",
   title = "{Domain decomposition methods-algorithms and theory}",
   series = "Springer Ser. Comput. Math.",
   volume = 34,
   publisher={Springer},
   year = 2005
   }

@article {MR3815657,
    AUTHOR = {Calvo, Juan G.},
     TITLE = {On the approximation of a virtual coarse space for domain
              decomposition methods in two dimensions},
   JOURNAL = {Math. Models Methods Appl. Sci.},
  FJOURNAL = {Mathematical Models and Methods in Applied Sciences},
    VOLUME = {28},
      YEAR = {2018},
    NUMBER = {7},
     PAGES = {1267--1289},
      ISSN = {0218-2025},
   MRCLASS = {65N30 (35J25 65N55)},
  MRNUMBER = {3815657},
       DOI = {10.1142/S0218202518500343},
       URL = {https://doi.org/10.1142/S0218202518500343},
}

@article {MR3913656,
    AUTHOR = {Calvo, Juan G.},
     TITLE = {An overlapping {S}chwarz method for virtual element
              discretizations in two dimensions},
   JOURNAL = {Comput. Math. Appl.},
  FJOURNAL = {Computers \& Mathematics with Applications. An International
              Journal},
    VOLUME = {77},
      YEAR = {2019},
    NUMBER = {4},
     PAGES = {1163--1177},
      ISSN = {0898-1221},
   MRCLASS = {65N55 (35J25 74S05)},
  MRNUMBER = {3913656},
MRREVIEWER = {Anh-Khoa Vo},
       DOI = {10.1016/j.camwa.2018.10.043},
       URL = {https://doi.org/10.1016/j.camwa.2018.10.043},
}

@article{NadeemJimack,
    author = {Nadeem, S. A. and Jimack, P. K.},
    title = "{Parallel implementation of an optimal two level additive Schwarz preconditioner for the 3-D finite element solution of elliptic partial differential equations}",
    journal = {Int. J. Numer. Methods Fluids},
    volume = {40},
    number = {12},
    pages = {1571-1579},
    doi = {https://doi.org/10.1002/fld.413},
    url = {https://onlinelibrary.wiley.com/doi/abs/10.1002/fld.413},
    eprint = {https://onlinelibrary.wiley.com/doi/pdf/10.1002/fld.413},
    year = {2002}
}

@article{UReading,
    author = {Schoeberl, Joachim and Melenk, Jens and Pechstein, Clemens and Zaglmayr, Sabine},
    year = {2008},
    month = {01},
    pages = {1-24},
    title = "{Additive Schwarz preconditioning for p-version triangular and tetrahedral finite elements}",
    volume = {28},
    journal = "IMA J. Numer. Anal.",
    doi = {10.1093/imanum/drl046}
}

@article{VEMcomp,
    author = {Frittelli, Massimo and Madzvamuse, Anotida and Sgura, Ivonne},
    year = {2024},
    month = {08},
    pages = {1393-1428},
    title = "{VEMcomp: a Virtual Elements MATLAB package for bulk-surface PDEs in 2D and 3D}",
    volume = {99},
    JOURNAL = {Numer. Algorithms},
    doi = {10.1007/s11075-024-01919-4}
}

@misc{mVEM,
      title="{mVEM: a MATLAB software package for the virtual element method}", 
      author={Yue Yu},
      year={2022},
      eprint={2204.01339},
      archivePrefix={arXiv},
      primaryClass={math.NA},
      url={https://arxiv.org/abs/2204.01339}, 
}

@article{BeiraoStab,
author = {Beir\~{a}o da Veiga, Louren\c{c}o and Lovadina, Carlo and Russo, Alessandro},
title = {Stability analysis for the virtual element method},
JOURNAL = {Math. Models Methods Appl. Sci.},
volume = {27},
number = {13},
pages = {2557-2594},
year = {2017},
doi = {10.1142/S021820251750052X},
URL = {https://doi.org/10.1142/S021820251750052X},
eprint = {https://doi.org/10.1142/S021820251750052X}
}

\end{document}